# Stable and Renewable: Assessing the Reliability of a Fully Renewable European Energy System

D. Franzmann[1,2,], N. Ludwig[1,3], J. Linßen[1], D. Stolten[2], H. Heinrichs[1,3*]

1 Forschungszentrum Jülich GmbH, Institute of Energy and Climate Research – Jülich Systems Analysis (ICE-2), 52425 Jülich, Germany

2 RWTH Aachen University, Chair for Fuel Cells, Faculty of Mechanical Engineering, 52062 Aachen, Germany

3 University of Siegen, Chair for Energy Systems Analysis, Department of Mechanical Engineering, 57076 Siegen, Germany

* corresponding author: h.heinrichs@fz-juelich.de

## Abstract

The transformation of the energy system has raised concerns about the reliability of fully renewable energy systems. We address this question for a 2050 European energy system using an economically optimal adequacy assessment. Our results show that a cost-optimal, fully renewable European system can be as reliable as a fossil-based one, with an average loss of load of only 0.03% due to variability in renewable generation. Outages primarily affect industrial and service sectors, while household supply remains largely uninterrupted. Regional differences in supply security emerge, with outages concentrated in countries with a low Value of Lost Load (VoLL). We demonstrate that system reliability can be fully ensured at negligible additional cost (+0.17%) by modestly increasing hydrogen turbine (+10%) and battery capacities (+15%) beyond the cost-optimal levels. We conclude that well-designed renewable energy systems are stable, with hydrogen-based backup being a key enabler of reliability.

## Keywords

energy security, partial outage, loss of load expectation, reliability standard, multi-day wind lulls, stability, adequacy assessment

# 1 Introduction

The transformation of energy systems towards high shares of variable renewable energy sources (RES) has raised the question of whether renewable energy systems can be reliable in the future. While conventional power plants adjust their generation to meet electricity demand, renewable energy systems rely on flexibility to match variable generation with demand. As a result, ensuring adequacy in renewable-dominated systems requires more sophisticated approaches to planning and investment decisions.

The variability of renewable energy sources is well understood. Wind power can experience widespread lulls [1], [2] for up to seven days [3], [4]. PV lulls last only for a maximum of two days, occur less frequently than wind lulls [5], and are also more predictable by diurnal and seasonal patterns [6]. Hydroelectric generation is mostly seasonal and varies greatly by region [3]. In contrast, the role of flexibility options has primarily been studied for specific technologies. For example, Weber et al. [5] and Boston et al. [7] show that ideally, a combination of backup power plants and storage facilities should be used to reduce outages at low cost. Boston et al. [7], Weber et al. [5], and Ryberg [8] indicate that outages can be avoided with sufficiently large storage and backup capacities, but economic limits exist to the expansion of these technologies. Grid expansion generally reduces the occurrence of outages [5], [7]. However, grid expansion cannot prevent all outages, either in a Europe without grid restrictions [8] or in a real-world scenario for Germany [9]. Overall, existing research demonstrates that technical solutions exist to balance variable generation [10]. However, the economic viability of such solutions remains at least in parts unresolved in the literature. Moreover, most studies rely on fixed assumptions about the deployment of renewable technologies, thereby implicitly defining the system's security of supply. This raises the question of how supply security emerges as an endogenous outcome of system design and cost optimization.

To be economically feasible, an economically optimal adequacy assessment must be considered. This approach allows the identification of economically optimal investment levels by balancing the cost of supply shortfalls, formally described as the Value of Lost Load (VoLL), against the cost of additional capacity, called the Cost of New Entry (CONE). When an economically optimal adequacy assessment is considered, the resulting energy system is cost-optimal in terms of energy supply and outage costs. This results in an economically optimal level of reliability, expressed as the optimal Loss of Load Expectation ($LOLE_{opt}$), defined in the following equation:

$$\text{LOLE}_{\text{opt}} = \frac{CONE}{VoLL}$$

This approach has been used by the EU since 2019 and is called the VOLL/CONE/RS approach, where RS stands for the resulting reliability standard [10]. The $LOLE_{opt}$ is calculated per bidding zone [11] and is used as an aim for the adequacy assessment and further measures for stabilization, like capacity markets [12]. In 2021, 10 of 27 EU states have activated these additional measures [13].

In this work, we answer the question of how stable a renewable energy system in the future can be if it is based only on renewable energy. Therefore, we use the economically optimal adequacy assessment in combination with a renewable energy system model to determine what outages will look like in the future.

## 2  Methodology

To answer the research question in the context of future renewable energy systems, we use an economically optimal adequacy assessment in combination with a high-resolution energy system model. We model the Cost of New Entry (CONE) based on the renewable energy system described in the following section, where free capacity expansion is allowed for optimal adequacy assessment. The value of lost load is derived based on the official EU values in Section 2.2.

### 2.1 Energy system design

In this paper, we model new entry costs using the high-resolution energy system model ETHOS.modelBuilder for renewable energy. It is based on the approach used in Franzmann et al. [14], which models a greenfield renewable energy system with sector coupling, including hydrogen back-up capacity and dispatchable renewable energy sources. We expanded it by a loss of load component based on the Value of Lost Load derived in the next section in more detail. The components considered for the energy system design are shown in Figure 1. The model represents the whole of Europe as an interconnected energy system based on the exchange of electricity and hydrogen. It consists of 250 regions based on the subnational GID-1, referring to the first-level administrative subdivisions, e.g., Federal States in Germany [15], aggregation with hourly resolution. The scope of the model is a renewable energy system in 2050. Therefore, the primary energy supply is limited to electricity based on ERA5 weather data [16]. For this study, the following flexibility options are included to allow the power system to avoid load losses: (1) generation flexibility through free curtailment and the use of dispatchable renewables, (2) temporal flexibility through Li-ion batteries, (3) regional flexibility based on electricity grids, and (4) sector coupling of hydrogen, which further allows temporal flexibility through long-term storage, regional flexibility based on hydrogen transport through pipelines, and production flexibility based on backup gas turbines. Therefore, our model includes all technical and operational flexibility options defined by Heider et al. [17], [18]. Openly available and spatially resolved data on the demand side flexibility, which do not exist by now, would be needed on a European scale to be included in future studies.

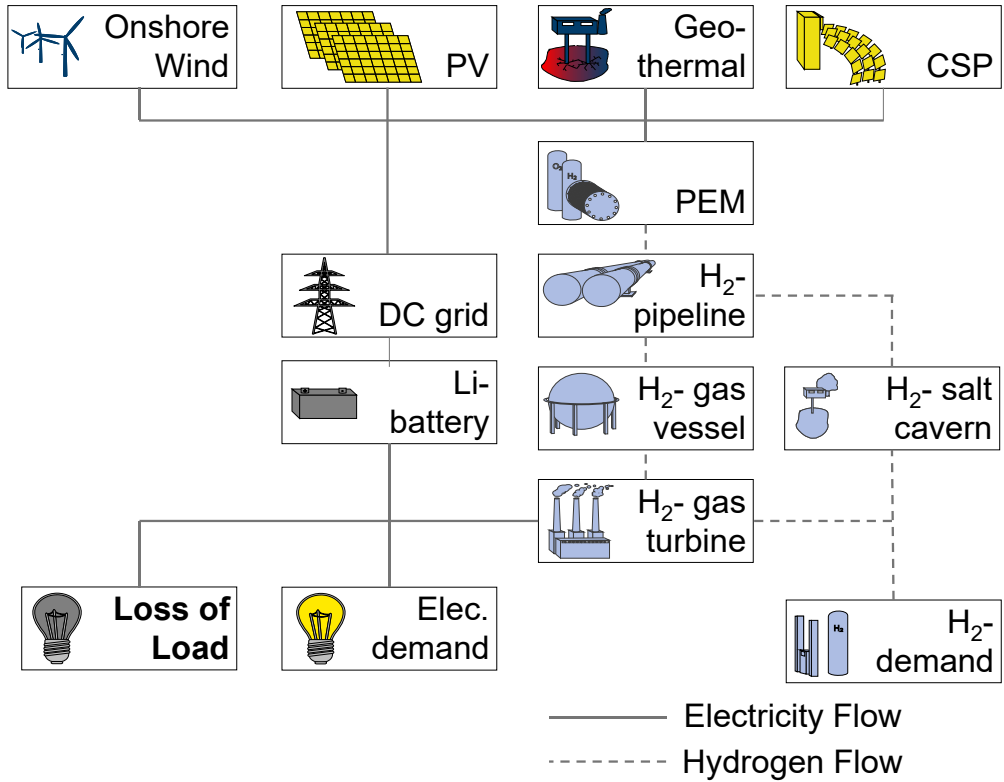

**Figure 1: Considered components within the energy system design approach**.

All capacities can expand endogenously in the model, constrained only by technical and geographical limitations, to address the full scope of the VoLL/CONE/RS [10] adequacy assessment. As a result, the model performs a cost-optimal adequacy assessment, resulting in an economically efficient level of supply security. We have modeled the economic trade-off between the supply cost $C_{System}$ and the loss of load (LOL) cost $C_{LOL}$ by adding the latter to the general cost minimization problem of techno-economic models [19]:

$$\min_{x} C_{System} + C_{LOL}$$
$$s.t.: Ax = b$$
$$Dx \leq e$$
$$C_{System} = c^T x$$
$$C_{LOL} = VoLL^T x$$

In this formulation, the loss of load expectation $\text{LOLE}_{\text{opt}}$ is included in the optimal solution $x^*$ of the free optimization variables $x$. $Ax = b$ and $Dx \leq e$ describe the equations and inequations from the physical formulation of the optimization problem.

## 2.2 Derive European VOLL for 2050

To model the VoLL/CONE/RS methodology, we used the official value of lost load from the EU estimated by CEPA in 2018 [20] and projected it to 2050 for model consistency. The projection is based on the national electricity demand $E_{el}$ and the gross domestic product GDP from the GCAM Net Zero 2050 scenario [21], [22], which was also used for the electricity demand for model consistency:

$$VoLL_{2050} = VoLL_{2020} \frac{GDP_{2050}/GDP_{2020}}{E_{el,2050}/E_{el,2020}}$$

This projection uses a common method to determine the VoLL based on a country's gross value added $GVA$ and electricity demand $E_{el}$: $VoLL = GVA/E_{el}$ [23]. CEPA only provides VoLL values for the 27 EU member states. Missing VoLL values were calculated directly from the country-specific GVA [24] and electricity demand [25] for 2020 and projected in the same way. The VoLL values for 2020 and 2050 are presented in Figure 2.

In addition to national differences, there are also large sectoral differences in VoLL. We model sectoral differences based on the available values from the Netherlands [26], Ireland [27], Germany [28], and Spain [29]. These are used to derive an average per sector: agriculture (22.1 EUR/kWh), services (4.1 EUR/kWh), households (19.9 EUR/kWh), industry (4.3 EUR/kWh), and transport (11.0 EUR/kWh) for an average VoLL of 7.3 EUR/kWh. Based on this distribution and the sectoral electricity demand from Franzmann et al. [14], the sectoral values are distributed to maintain the average country VoLL from Figure 2.

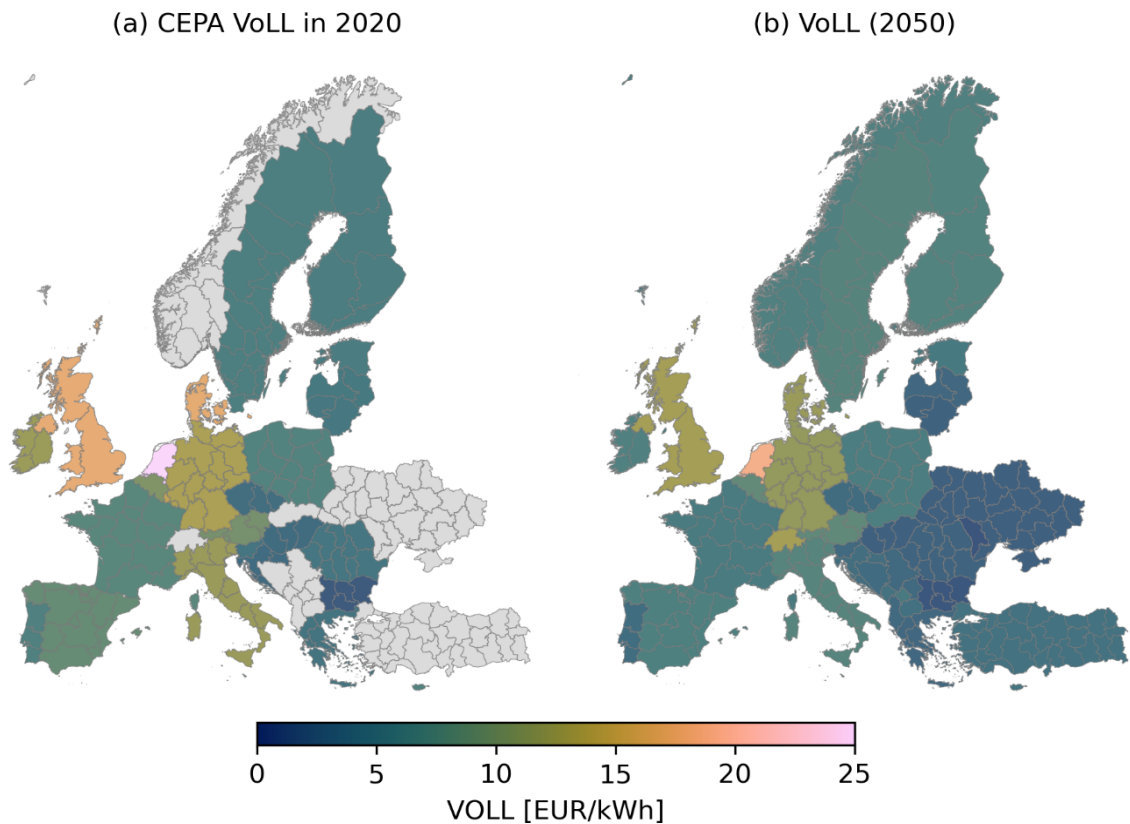

**Figure 2:** (a) Official CEPA VoLL value for the 27 EU Member States where officially calculated and (b) projected average VoLL value for 2050 based on GDP and electricity demand. The average VoLL in Europe is expected to decrease as electricity demand increases more than GDP due to the electrification of the energy system based on the underlying Net Zero 2050 scenario [21], [22].

## 3  Results

Security of supply depends heavily on the design of the energy system. Therefore, we first describe what the economically optimal energy system design looks like. Then, we explore what outages might look like in a future renewable Europe based on the optimal adequacy assessment. Finally, we derive properties of stable renewable energy systems.

### 3.1 Energy system design

The energy system derived for a 100% renewable Europe in 2050 consists mainly of wind generation (60% of generation) and PV generation (39%). Most of the solar generation takes place in Spain, Greece, and southern Italy, while the wind generation is distributed over the northern UK, eastern Europe, Scandinavia, and Spain (see Figure 3). Interestingly, these high-generation regions also have the highest capacity for PEM operation, while hydrogen gas turbines for re-electrification are located in the demand centers of Central Europe. In total, 2709 GW of PV, 2585 GW of onshore wind, 1116 GW of PEM electrolysis, and 501 GW of back-up hydrogen gas plants will be built. In contrast, the EU-27 and Europe have a capacity of 269 GW of solar [30] and 285 GW of wind [31] installed in 2023 and 2024, respectively. Both correspond to about 10% of the expansion targets in this study. The total annual cost of the energy system is 91.6 EUR/MWh, of which 44% is for onshore wind, 14% for PV, 22% for hydrogen infrastructure, including PEM electrolysis, hydrogen storage, pipelines, and $H_2$ gas turbines. The electricity grid accounts for 15% of the total annual cost, whereas Li-ion batteries contribute only 1% of the total cost, as most of the storage (98%) is handled by hydrogen.

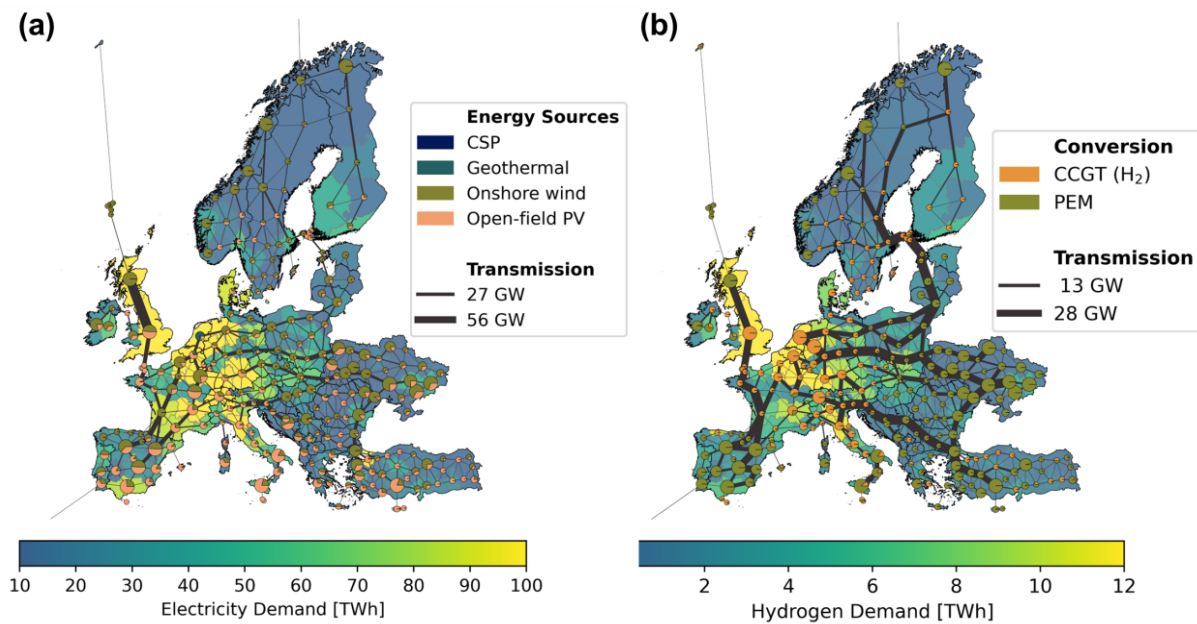

**Figure 3: Energy system capacities and demand of a fully renewable Europe in 2050 with optimal adequacy assessment. (a) shows the electricity sector, while (b) shows the hydrogen sector.**

## 3.1 Economic efficient outages of renewable energy systems in Europe

In a renewable Europe in 2050, based on a cost-optimal adequacy assessment, the distribution of load losses is highly variable. On average, 0.03% of the European load is lost, with a range from 0 to 0.54% (see Figure 4a). The lowest load loss, below 0.04%, is observed in Western and Northern Europe and Turkey. The only regions in these parts of Europe with noticeable load loss are the edges of Portugal and northern Norway, below 0.14%. High load losses occur in Eastern Europe, with load losses of up to 0.54% in Bulgaria and Hungary, 0.4% in Moldova, and 0.36% in Ukraine. Comparing the loss of load with the value of lost load from Figure 2, the high impact of the value of lost load on the loss of load can be seen. In countries with a VoLL below 3 EUR/kWh, the highest losses are observed, with up to 0.5% of energy not delivered (Bulgaria, Hungary, Lithuania, Moldova, Romania, Ukraine). In the intermediate range of 3 to 5.5 EUR/kWh are the countries with medium load loss, such as most of Eastern Europe, as well as Greece and Portugal. Above 5.5 EUR/kWh are the most stable countries in Western Europe, with 80% of countries having load losses below < 0.001%. Therefore, the value of load loss has the greatest impact on the security of supply in a future Europe.

The duration of outages shows a homogeneous characteristic as presented in Figure 4b. Almost all partial outages last between 5 and 10 hours. On a European average, renewable energy systems with an economically optimal adequacy assessment experience 2.2 partial outages per region with a duration of 5 hours. Partial outages longer than 10 hours occur only once every 12 years per region on average in Europe. Regions in the lowest VoLL category experience more than 8 partial outages per region, while countries in the highest VoLL category experience between 0 and 1 outages per year.

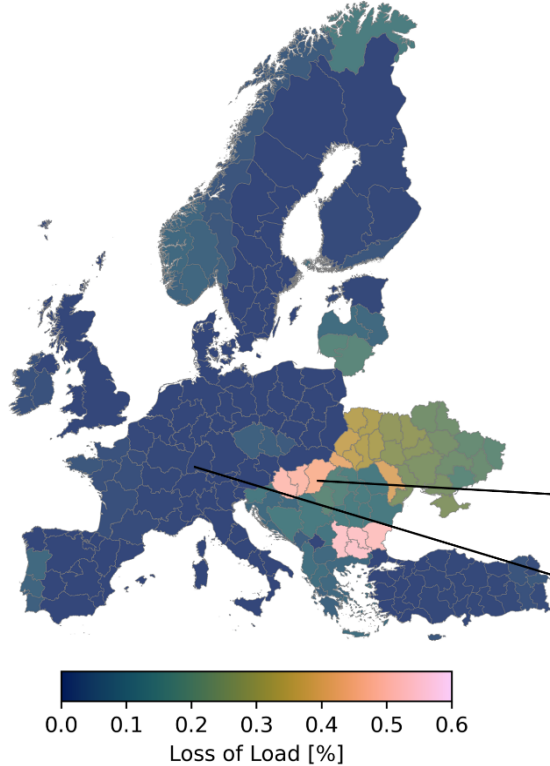
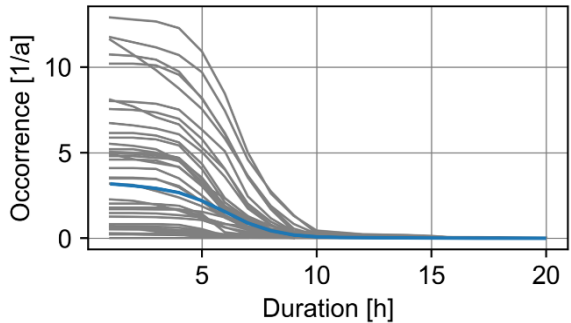
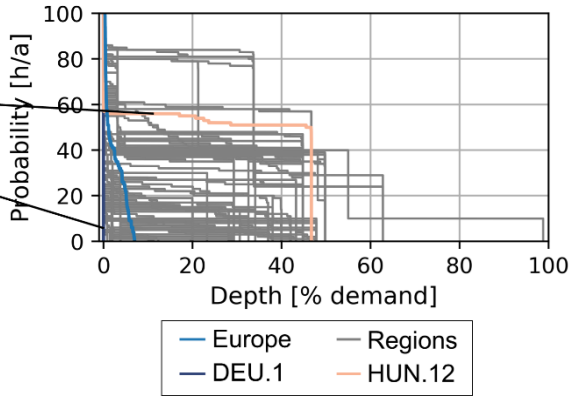

**Figure 4: Optimized energy losses in Europe over 20 different weather years. (a)** The average power loss per region is less than 0.6% per year, with the highest power losses in Eastern Europe, Norway, and Portugal. **(b)** The average number of outages per region per year shows an average of 7 outages per year, with a variation from 0.0 in Germany to 12.9 in Bulgaria. The frequency of outages decreases sharply between 5 and 10 hours. **(c)** Outages do not usually lead to complete blackouts. On average in Europe, only a small number of regions experience power cuts of more than 50% of regional demand.

Energy losses are low across Europe (see Figure 4c). The maximum instantaneous power not provided at the European level is 7%, and the European power loss of 5% occurs for 20 hours per year on a 20-year average. As described above, there are significant differences between Western and Eastern Europe. While Baden-Württemberg in Germany (DEU.1) has no losses at all, Hajdú-Bihar in Hungary (HUN.12) has 48% load loss for more than 50 hours per year. 58 out of 250 regions, mostly in Eastern Europe, experience outages of more than 40% of their load, but only 5 regions in Bulgaria and Moldova experience more than 50% of their load not supplied. In Western and Northern Europe, 150 out of 250 reach a maximum load loss of less than 2%.

In Western Europe, outages due to fluctuations in renewable energy are sparse, with maximum outages of less than 2%. In contrast, Eastern Europe experiences much deeper outages of 50% electricity consumption, and outages can occur more than 8 times per year. In renewable Europe, the duration of brownouts is usually less than 10 hours.

### 3.2 Renewable energy systems export lulls to low VoLL regions

Typically, a load-loss event occurs during a wind generation lull. The wind outage results in a high residual load, defined as the electricity demand that is not met by direct local renewable electricity generation. If the high residual load cannot be compensated for by the flexibility option, an outage occurs. Within Europe, we observe two different outage patterns depending on whether the wind lull occurs in Western or Eastern Europe (see Figure 5).

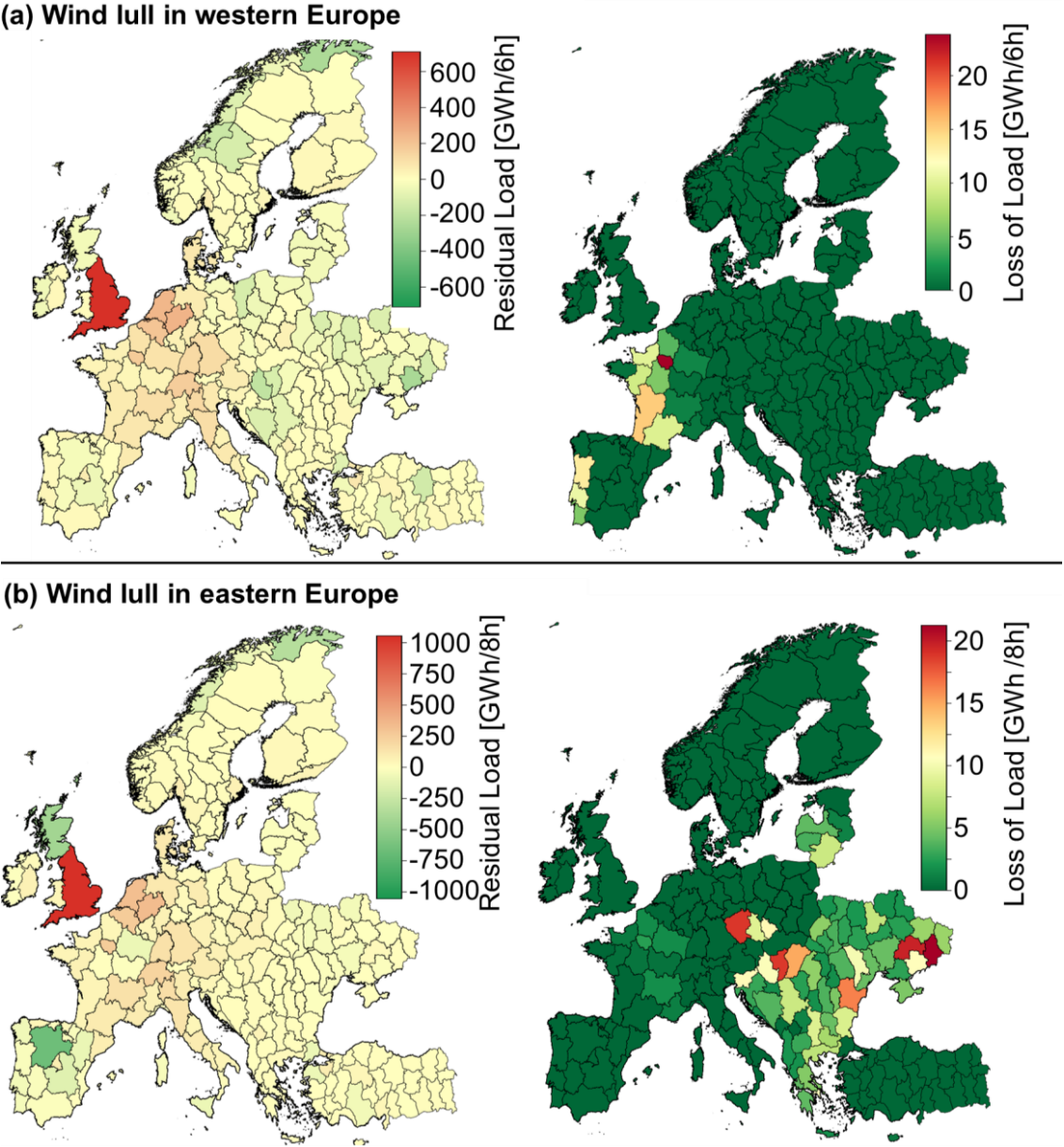

**Figure 5: The two different outage patterns in a fully renewable European energy system: (a) Spain experiences a lull in wind generation with reduced energy exports to demand centers in Central Europe. The partial outages are observed in Portugal and France. (b) Eastern Europe experiences a wind generation lull, with reduces its energy exports to the demand centers in Central Europe. The partial outages are**

**observed throughout Eastern Europe. (a) is observed during hours 8104 to 8110 of the 2018 weather data, and (b) is observed during hours 855 to 863 of the same weather year.**

This is due to the location of high VoLL demand centers in Europe. They are located on an axis from England to northern Italy and do not allow any eastward or westward energy transport perpendicular to this axis during times of energy scarcity, as their supply is prioritized due to their highest VoLLs. This observation underscores the complexity of large interconnected systems and indicates that even well-integrated networks can experience regionally differentiated stress modes. The two patterns are described below:

In the example for the Western European case, the large wind capacities of the Iberian Peninsula and northern Great Britain experience a lull in wind generation (see Figure 5a). Under normal operating conditions, regions such as Spain and northern England have significant renewable generation surpluses, with residual loads as low as -125 GW, and act as major exporters of electricity and hydrogen (see Figure 5b, where Western Europe's electricity generation is normal). Without Western Europe's energy exports, there is a deficit. Although England has the highest absolute deficit of over 600 GWh, it does not face a load loss due to its high VoLL of 13.27 EUR/kWh. On the other hand, Portugal and France suffer from load shedding despite comparatively smaller deficits. These blackouts are mainly due to lower VoLL values of 3.65 EUR/kWh in Portugal and 5.60 EUR/kWh in France, which make load shedding more cost-optimal than expensive dispatch alternatives. Although Eastern Europe sheds up to 47 TWh of electricity during the event, enough to cover the 1 TWh shortfall in Western Europe, grid congestion prevents effective redistribution. Instead, a significant portion of this surplus (1.7 TWh) is converted to hydrogen and shipped to Central Europe, where it remains unusable due to local hydrogen-to-electricity capacity constraints resulting from the Europe-wide cost optimization.

The second cause of energy system outages is wind power outages in Eastern Europe (see Figure 5b). The eight-hour outage with low wind generation in Eastern Europe results in a shortfall of 627 GWh in Europe. High residual load occurs mainly in the high-energy-demand regions of Central Europe. Due to Eastern Europe's role as an exporter of wind power, there is no negative residual load in Eastern Europe despite the active wind lull. During the lull, the individual regions of Eastern Europe have an average electricity surplus of 1.9 GW. Nevertheless, partial outages occur precisely in the regions with an electricity surplus (see Figure 5b). For example, Ukraine has a surplus of 350 GWh of electricity during the eight hours considered, but exports 510 GWh of electricity. This results in a shortfall of 175 GWh in Ukraine, as 34% of the electricity is used for export from Ukraine.

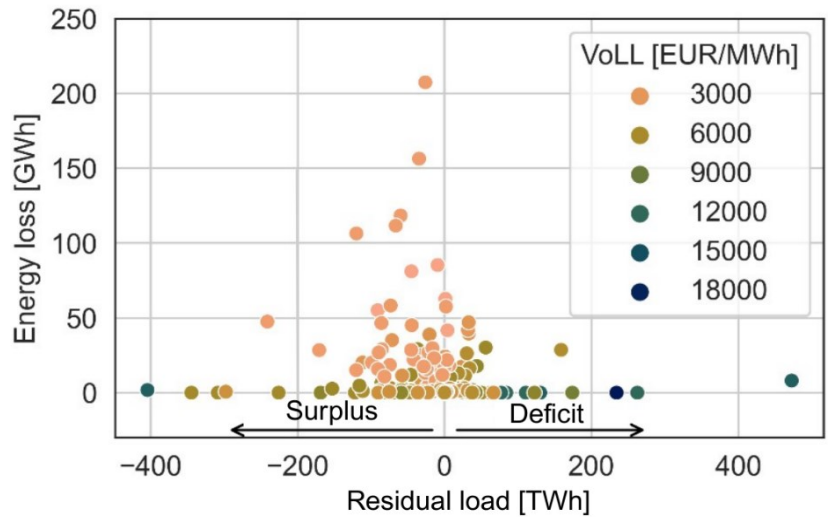

**Figure 6: The loss of load per region in Europe as a function of the value of lost load and the residual load shows that outages do not occur in regions with an energy deficit, but rather in regions with a low VoLL. Due to the correlation between high demand and high VoLL, outages tend to occur in regions with an abundant power supply.**

Both outage types, East Europe and West Europe, show a pattern of load loss being diverted to the region with the lowest VoLL during operation. This is shown in Figure 6, which shows the loss of load as a function of regional residual load and VoLL. Few outages occur in regions with an actual deficit. Instead, outages occur predominantly in regions with a VoLL of less than 4 EUR/kWh. Even additional costs due to energy losses for the conversion of electricity to hydrogen are accepted, which further increases the shortfall during power outages. This shows that there is a systematic export of losses to regions with a low VoLL. These regions are in Eastern Europe, Portugal, and France.

### 3.3 Techno-economic limitations leading to partial outages

This section describes the reasons for outages. The future European renewable energy system in this work is mainly powered by wind (60% of the energy) and PV (39%). The power outages are mainly caused by multi-day lulls in wind generation, as marked by the grey rectangle in Figure 7c. While solar generation can temporarily mitigate these shortfalls during the day, outages occur primarily in the evening and night due to the absence of solar irradiance (see Figure 7d). This diurnal pattern explains the prevalence of short-duration outages below 10h, which are often interrupted by daytime solar input. However, the root cause is often extended periods of low wind generation lasting up to 7 days.

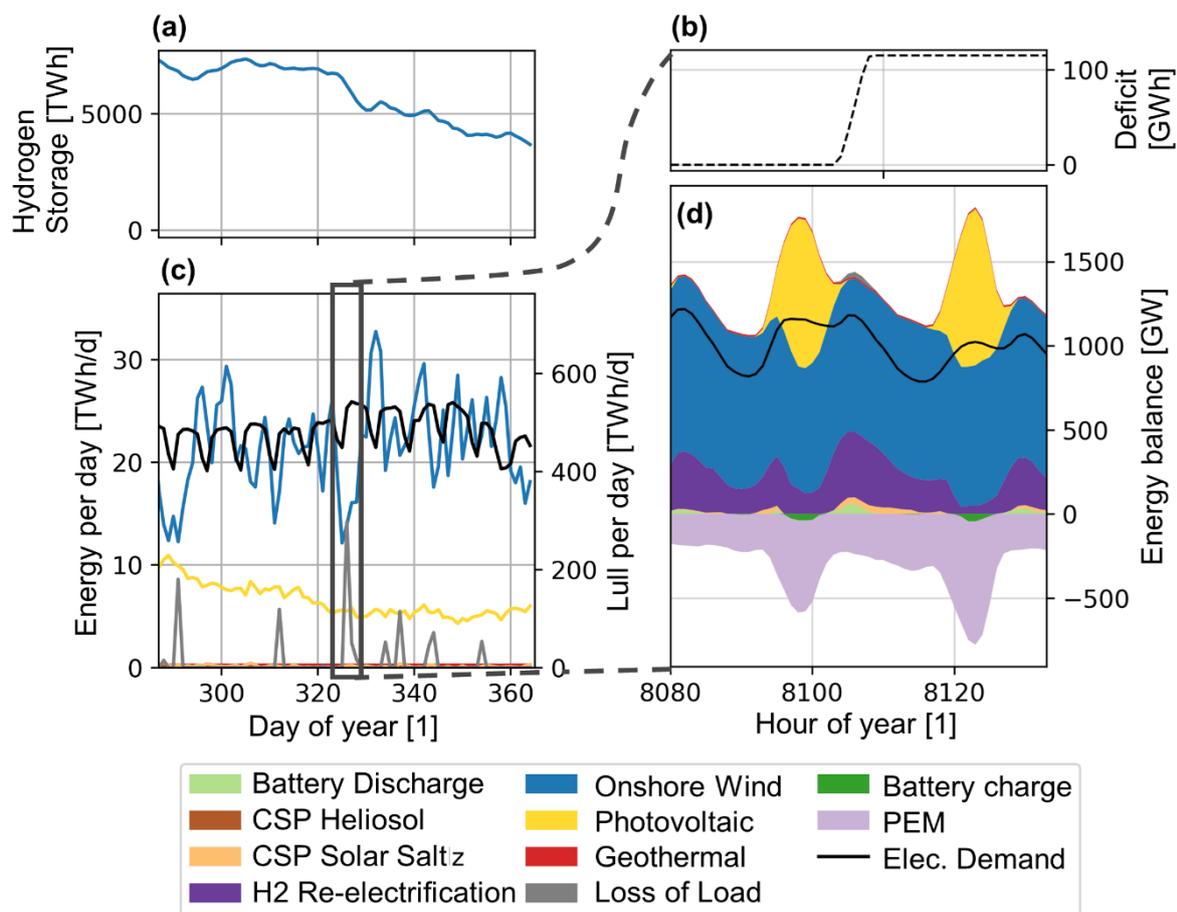

Figure 7: Daily and hourly time series of the European energy system during a theoretical outage event from 8080h-8136h, based on the weather year 2018. (c) shows that the reason for the outage is a generation lull of wind power, which cannot be compensated for and leads to a 100GWh deficit over 6 hours (b) in the right figure. (a) shows that the stored amount of hydrogen is not limited, but rather the available backup power (see d).

In the example in Figure 7, even though the wind electricity generation is reduced by -50%, 80% of the energy is still supplied by wind power during the wind lull. The energy system compensates for most of the reduced generation by storing energy. Energy is shifted from the daytime surplus to the evening deficit using hydrogen storage (98%) and batteries (2%). The backup hydrogen energy is not limited at any point because the hydrogen storage level is always above 5000 TWh during the lull (Figure 7a). Instead, the backup power of the reconversion and the battery capacity limit the operation. During the lull, all hydrogen gas turbines are in full operation, and the batteries are fully discharged. Therefore, the economic trade-off between the cost of security of supply and the cost of energy losses is directly manifested in the expansion of hydrogen gas turbines for backup power and, to 2%, in the expansion of batteries. This means that the security of supply of renewable energy systems can be directly set by the expansion of hydrogen gas turbines for backup power.

### 3.4 Fully stable renewable energy systems are achievable at low costs

Instead of accepting the security of supply indirectly set by the VoLL, stable renewable energy systems can be built. Here, we define a stable energy system by an outage of less than 1 hour per year and consumer, which reflects the temporal resolution of the model and thus represents the highest stability observable. A stable European renewable energy system can be achieved by 3 major steps (see Figure 8): (1) an additional 9.8 GW (+2%) of hydrogen gas turbines are added to the system to reduce the reconversion constraints described above, (2) 8 GW (-0.3%) of wind turbine capacity is replaced by 5.2 GW (+0.2%) of PV capacity to reduce the impact of wind lulls as the reason for outages, and (3) storage, especially batteries, is expanded by +15%. In total, the stable energy system can be achieved at no relevant higher cost (0.17% additional cost). Therefore, a stable, fully renewable energy system is achievable within technical and economic constraints.

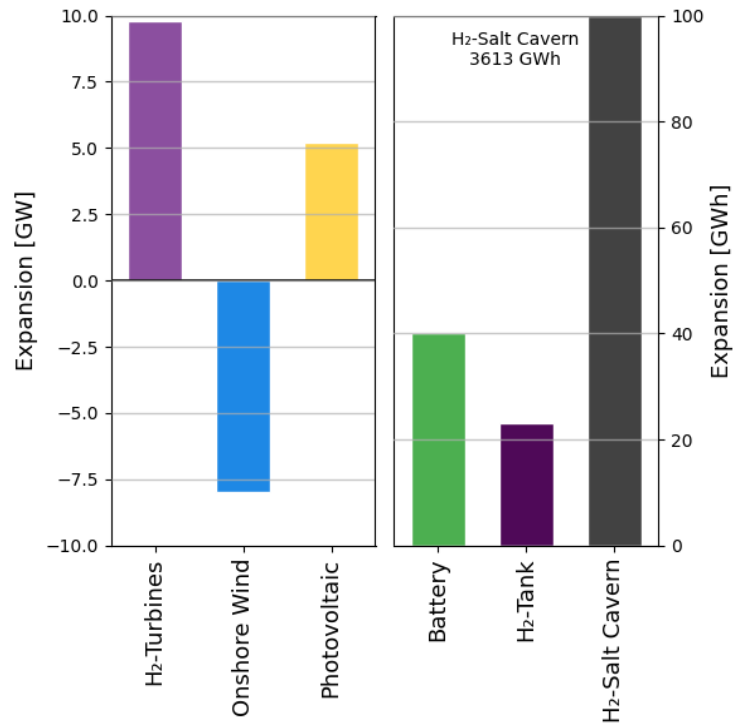

**Figure 8: Additional capacity of a stable power system compared to an economically optimal power system.**

## 4 Discussion

The results show that the European renewable energy system based on a value of lost load approach will be stable, with an average load loss probability of 0.03%, matching the current reliability standards of 2.77h/a (=0.03%) for Germany and Luxembourg [32] and below the current reliability standard of 8h/a (=0.10%) for Lithuania [33]. It should be noted that we have modeled a technically fully operational system, which does not include outages due to technical failures, which are included in the German reliability standard. On the other hand, the technical redundancy of renewable energy systems increases as generation becomes more decentralized and transmission capacity is expanded. Therefore, the reliability of a fully renewable European energy system may be slightly reduced in comparison to the current standard, but the difference is expected to remain minor. In Western Europe, the additional load losses from renewable energy systems are zero and do not add any additional losses to the power system.

The average power loss during outages reduces supply by less than half, meaning that at least 50% of demand is served. Therefore, not all sectors are equally affected by power outages, as some power is always supplied. In particular, the industrial and service sectors, which have the lowest sectoral VoLL, are affected by brownouts, while households usually do not have any outages at all. Therefore, these sectors may have to adapt to flexible energy use in the future, although these events only occur between 0 and 12 times per year, depending on the country, challenging the economic viability of measures for demand side management. In order to have a controlled brownout, technologies for early warning as well as load shedding of the sectors are needed to avoid unstable situations in the energy system.

We find that there are large regional differences in outages between Western and Eastern Europe and show that these are due to a systematic shift in exports from high-VoLL countries in Central Europe to low-VoLL countries in Eastern Europe. As the reliability of power supply is a factor for certain industries, this can lead to systematic disadvantages in industrial development and power supply for economically weaker countries with lower VoLL. Therefore, this study recommends embedding market mechanisms that limit the export of energy in times of energy scarcity to avoid shifting energy losses to low-VoLL regions.

Finally, this study also shows that adding backup power in the form of a hydrogen gas turbine and batteries can provide stable renewable power at almost no additional cost. This would further reduce the cost of managing brownouts and reduce the risk of systematically disadvantaging economically weak countries by exporting brownouts. It could also improve the acceptance of the energy transition by providing a high standard of energy supply.

Due to the integral consideration of supply interruptions, this study yields deviating results compared to other studies focusing on the European energy system, hindering direct comparison. This is further aggravated by deviating assumptions in terms of system boundaries, energy demands, renewable energy potentials, utilized weather years, as well as techno-economic parameters for which no consensus exists within the literature. This leads, for example, to higher wind and lower battery capacities in comparison to the study from Dunkel et al. [10], mainly driven by lower electricity demands due to hydrogen imports. However, neglecting hydrogen imports as a potential but uncertain flexibility option allows this study to investigate the ability of renewable energy systems in Europe to provide a stable supply even under challenging geopolitical conditions.

The limitations of this study are that CONE and VoLL in the future are not well known, but are based on projections. In particular, the CONE depends on a further reduction in the cost of renewable technologies, which must be technically feasible. Therefore, the results of this study are based on a further reduction of technology costs until 2050. In addition, this study allows for an expansion of renewable energy technologies that are only limited by the land availability for generation units. In particular, additional constraints on transport capacity could alter the results, although open data in this regard are sparse [34]. However, as the scope of this study is to investigate theoretical techno-economic limits, this can be seen as the next step.

## 5 Conclusions

This study shows that a fully renewable European energy system by 2050 is both technically and economically feasible, achieving a high level of reliability with an average load loss of only 0.03%. The resulting outages are always shorter than 10 hours. While Eastern Europe has more than 8 partial outages per region per year, Western Europe has only between 0 and 1 outage per year and region. The outages are mostly partial. In the regions of Western Europe, 150 out of 250 regions have a maximum power outage of only 2%. Eastern Europe has higher outages of over 40% due to its lower VoLL. Countries with a VoLL below 3.3 EUR/kWh, mostly in Eastern Europe, have higher load-losses than countries with a VoLL above 5.5 EUR/kWh in Western Europe. Importantly, the VoLL-based adequacy assessment leads to a systematic redistribution of energy shortages towards economically weaker regions, raising concerns about reliability equality. Outages are primarily driven by multi-day wind lulls and are limited by hydrogen gas turbine back-up capacity and battery storage capacity rather than hydrogen back-up energy.

Finally, we find that renewable energy systems can be built at an incremental cost of +0.17% by further expanding hydrogen gas turbines and battery capacity. These two technologies are important in determining the level of security of renewable energy systems and can be adapted retroactively to the expansion of renewable energy systems. We therefore recommend that policymakers and energy system planners aim directly at building stable energy systems and use hydrogen gas turbines and battery storage to fine-tune the desired level of security after the massive expansion of renewable energy that is needed.

## Acknowledgements

This work was mainly supported by the Helmholtz Association under the program "Energy System Design". This work was partly funded by the European Union (ERC, MATERIALIZE, 101076649). Views and opinions expressed are, however, those of the authors only and do not necessarily reflect those of the European Union or the European Research Council Executive Agency. Neither the European Union nor the granting authority can be held responsible for them.

## Declaration of generative AI and AI-assisted technologies in the writing process

During the preparation of this work, the authors used ChatGPT (GPT-4) and DeepL Write in order to improve clarity and language. After using this tool/service, the authors reviewed and edited the content as needed and take full responsibility for the content of the published article.

## Authors contributions

The following abbreviations are used: Conceptualization (C), Methodology (M), Formal Analysis (FA), Writing - original draft (WO), Writing - review and editing (WR), Supervision (S), Project administration (PA), Funding acquisition (FU), Visualization (V).

DR: C, M, FA, WO, WR, PA, V;

NL: WO, WR, V;

JL: FU;

DS: S;

HH: C, M, WO, WR, PA, FU, S

All authors have read and agreed to the published version of the manuscript.

## References

[1] K. Z. Rinaldi, J. A. Dowling, T. H. Ruggles, K. Caldeira, and N. S. Lewis, "Wind and Solar Resource Droughts in California Highlight the Benefits of Long-Term Storage and Integration with the Western Interconnect," *Environ. Sci. Technol.*, vol. 55, no. 9, pp. 6214–6226, 2021, doi: 10.1021/acs.est.0c07848.

[2] N. Ohlendorf and W.-P. Schill, "Frequency and duration of low-wind-power events in Germany," *Environ. Res. Lett.*, vol. 15, no. 8, p. 084045, 2020, doi: 10.1088/1748-9326/ab91e9.

[3] D. Raynaud, B. Hingray, B. François, and J. D. Creutin, "Energy droughts from variable renewable energy sources in European climates," *Renew. Energy*, vol. 125, pp. 578–589, 2018, doi: 10.1016/j.renene.2018.02.130.

[4] N. Otero, O. Martius, S. Allen, H. Bloomfield, and B. Schaefli, "A copula-based assessment of renewable energy droughts across Europe," *Renew. Energy*, vol. 201, pp. 667–677, 2022, doi: 10.1016/j.renene.2022.10.091.

[5] J. Weber, J. Wohland, M. Reyers, J. Moemken, C. Hoppe, J. G. Pinto, and D. Witthaut, "Impact of climate change on backup energy and storage needs in wind-dominated power systems in Europe," *PLoS One*, vol. 13, no. 8, p. e0201457, 2018, doi: 10.1371/journal.pone.0201457.

[6] K. Engeland, M. Borga, J.-D. Creutin, B. François, M.-H. Ramos, and J.-P. Vidal, "Space-time variability of climate variables and intermittent renewable electricity production – A review," *Renew. Sustain. Energy Rev.*, vol. 79, pp. 600–617, 2017, doi: 10.1016/j.rser.2017.05.046.

[7] A. Boston, G. D. Bongers, and N. Bongers, "Characterisation and mitigation of renewable droughts in the Australian National Electricity Market," *Environ. Res. Commun.*, vol. 4, no. 3, p. 031001, 2022, doi: 10.1088/2515-7620/ac5677.

[8] D. S. Ryberg, *Kalte Dunkelflauten der zukünftigen Potentiale der Wind- und Solarenergie in Europa*. RWTH Aachen University, 2020. doi: 10.18154/RWTH-2020-10968.

[9] F. Huneke, C. Linkenheil, and M. Niggemeier, *KALTE DUNKELFLAUTE: ROBUSTHEIT DES STROMSYSTEMS BEI EXTREMWETTER*. Greenpeace Energy eG, 2017.

[10] Philipp Dunkel, Theresa Klütz, Jochen Linßen, and Detlef Stolten, "Towards Hydrogen Autarky? Evaluating Import Costs and Domestic Competitiveness in European Energy Strategies." [Online]. Available: https://arxiv.org/abs/2510.04669

[11] Europäischen Union, "VERORDNUNG (EU) 2019/943 DES EUROPÄISCHEN PARLAMENTS UND DES RATES vom 5. Juni 2019 über den Elektrizitätsbinnenmarkt." 2019. [Online]. Available: http://data.europa.eu/eli/reg/2019/943/oj

[12] ACER, "Methodology for calculating the value of lost load, the cost of new entry and the reliability standard." 2022. [Online]. Available: https://acer.europa.eu/sites/default/files/documents/Decisions_annex/ACER%20Decision%2023-2020%20on%20VOLL%20CONE%20RS%20-%20Annex%20I.pdf

[13] European Comission, "Guidance for Member States on implementation plans pursuant Art. 20 (3)-(8) of Regulation (EU) 2019/943 ('Market Reform Plans')." 2020. [Online]. Available: https://energy.ec.europa.eu/document/download/674cbf1d-9bc5-4ccb-b634-a57f300eb900_en?filename=market_reform_plan_guidance

[14] European Union Agency for the Cooperation of Energy Regulators, "Security of EU electricity supply in 2021: Report on Member States approaches to assess and ensure adequacy." 2022. [Online]. Available: https://acer.europa.eu/Publications/ACER_Security_of_EU_Electricity_Supply_2021.pdf

[15] D. Franzmann, Christoph. Winkler, P. Dunkel, M. Stargardt, A. Burdack, S. Ishmam, J. Linssen, D. Stolten, and H. Heinrichs, "Impact of Spatial Aggregation on Global Renewable Energy System with ETHOS.modelBuilder," 2025, [Online]. Available: https://papers.ssrn.com/sol3/papers.cfm?abstract_id=5252316

[16] "Database of Global Administrative Areas, Version 4.1," GADM (or "Global Administrative Areas"), GADM, 2024. [Online]. Available: https://gadm.org/

[17] C3S, "ERA5 hourly data on single levels from 1940 to present." Copernicus Climate Change Service (C3S) Climate Data Store (CDS), 2018. doi: 10.24381/CDS.ADBB2D47.

[18] P. D. Lund, J. Lindgren, J. Mikkola, and J. Salpakari, "Review of energy system flexibility measures to enable high levels of variable renewable electricity," *Renew. Sustain. Energy Rev.*, vol. 45, pp. 785–807, 2015, doi: 10.1016/j.rser.2015.01.057.
[19] A. Heider, R. Reibsch, P. Blechinger, A. Linke, and G. Hug, "Flexibility options and their representation in open energy modelling tools," *Energy Strategy Rev.*, vol. 38, p. 100737, 2021, doi: 10.1016/j.esr.2021.100737.
[20] T. Klütz, K. Knosala, J. Behrens, R. Maier, M. Hoffmann, N. Pflugradt, and D. Stolten, "ETHOS.FINE: A Framework for Integrated Energy System Assessment," *J. Open Source Softw.*, vol. 10, no. 105, p. 6274, Jan. 2025, doi: 10.21105/joss.06274.
[21] Cambridge Economic Policy Associates Ltd, *STUDY ON THE ESTIMATION OF THE VALUE OF LOST LOAD OF ELECTRICITY SUPPLY IN EUROPE*. AGENCY FOR THE COOPERATION OF ENERGY REGULATORS, 2018. [Online]. Available: https://www.acer.europa.eu/en/Electricity/Infrastructure_and_network%20development/Infrastructure/Documents/CEPA%20study%20on%20the%20Value%20of%20Lost%20Load%20in%20the%20electricity%20supply.pdf
[22] K. Calvin, P. Patel, L. Clarke, G. Asrar, B. Bond-Lamberty, R. Y. Cui, A. Di Vittorio, K. Dorheim, J. Edmonds, C. Hartin, *et al.*, "GCAM v5.1: representing the linkages between energy, water, land, climate, and economic systems," *Geosci. Model Dev.*, vol. 12, no. 2, pp. 677–698, 2019, doi: 10.5194/gmd-12-677-2019.
[23] NGFS, "NGFS Climate Scenarios for central banks and supervisors." 2021. [Online]. Available: https://www.ngfs.net/sites/default/files/media/2021/08/27/ngfs_climate_scenarios_phase2_june2021.pdf
[24] M. J. Sullivan and D. M. Keane, *Outage cost estimation guidebook*. U. S. Department of Energy, 1995. [Online]. Available: https://www.osti.gov/biblio/239294
[25] WorldBank, "GDP (current US$)." Economic Policy & Debt: National accounts: US$ at current prices: Aggregate indicators, Sep. 13, 2023.
[26] U.S. Energy Information Administration, "Electricity generation, capacity, and sales in the United States." U.S. Energy Information Administration, 2022. [Online]. Available: https://www.eia.gov/energyexplained/electricity/electricity-in-the-us-generation-capacity-and-sales.php
[27] M. de Nooij, C. Koopmans, and C. Bijvoet, "The value of supply security," *Energy Econ.*, vol. 29, no. 2, pp. 277–295, 2007, doi: 10.1016/j.eneco.2006.05.022.
[28] E. Leahy and R. S. J. Tol, "An estimate of the value of lost load for Ireland," *Energy Policy*, vol. 39, no. 3, pp. 1514–1520, 2011, doi: 10.1016/j.enpol.2010.12.025.
[29] A. J. Praktiknjo, "Stated preferences based estimation of power interruption costs in private households: An example from Germany," *Energy*, vol. 76, pp. 82–90, 2014, doi: 10.1016/j.energy.2014.03.089.
[30] P. Linares and L. Rey, "The costs of electricity interruptions in Spain. Are we sending the right signals?," *Energy Policy*, vol. 61, pp. 751–760, 2013, doi: 10.1016/j.enpol.2013.05.083.
[31] "The Rise of Solar PV in the EU - key facts."
[32] Wind Europe, "Wind Europe. Wind energy in Europe: 2024 Statistics and the outlook for 2025-2030 2025." [Online]. Available: https://www.rinnovabili.it/wp-content/uploads/2025/02/Wind-energy-in-Europe-2024.pdf
[33] Bundesministerium für Wirtschaft und Klimaschutz, "Vorschlag der Regulierungsbehörden zum Zuverlässigkeitsstandard." [Online]. Available: https://www.bundeswirtschaftsministerium.de/Redaktion/DE/Downloads/V/vorschlag-der-regulierungsbehoerden-zum-zuverlaessigkeitsstandard.pdf?__blob=publicationFile
[34] ACER European Union Agency for the Cooperation of Energy Regulators, "Security of EU electricity supply in 2021: Report on Member States approaches to assess and ensure adequacy," 2022. [Online]. Available: the Cooperation of Energy Regulators. SeReprt on Member States approaches thttps://acer.europa.eu/sites/default/files/documents/Publications/ACER_Security_of_EU_Electricity_Supply_2021.pdf.

[35] W. Medjroubi, U. P. Müller, M. Scharf, C. Matke, and D. Kleinhans, "Open Data in Power Grid Modelling: New Approaches Towards Transparent Grid Models," *Energy Rep.*, vol. 3, pp. 14–21, Nov. 2017, doi: 10.1016/j.egyr.2016.12.001.

# Appendix

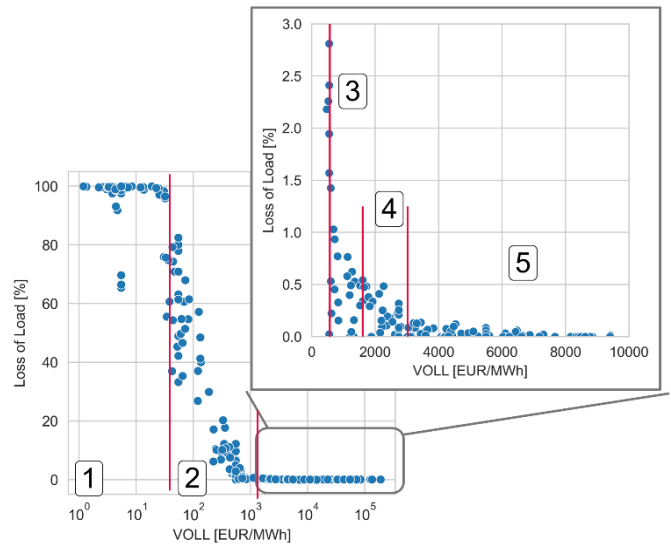

**Figure 9: Loss of Load (LOL) as a function of the Value of Lost Load (VoLL). Five characteristic regions are distinguished according to VoLL, with higher VoLL values leading to lower LOL.**

To quantify the effect of the Value of Lost Load (VoLL) on electricity supply stability, VoLL is varied from 0.001 to 10 times the reference value. As shown in Figure 9, this variation leads to a wide spectrum of outcomes—from complete system collapse to uninterrupted electricity supply—highlighting VOLL as a key parameter in system adequacy assessments. For interpretability, the results are divided into five characteristic VOLL intervals:

**(1) Region I (VOLL < 0.025 EUR/kWh):** At such low VOLL values, the model finds it more economical to accept complete electricity supply interruptions rather than invest in any generation or grid infrastructure. No European countries fall within this range.

**(2) Region II (0.025–0.5 EUR/kWh):** This interval marks a transition from full outages to nearly full supply. As VOLL increases above the levelized costs of renewables, investment in infrastructure becomes economically viable. However, no European countries in the baseline scenario exhibit VOLL values within this range.

**(3) Region III (0.5–3 EUR/kWh):** This range includes regions with partially stabilized systems but still significant outages. In this range, the median loss of load is at 0.36%, which is about 32 hours per year. From the base VoLL values, six Eastern European countries fall under this scenario and show a loss of loads of up to 0.5% (46h of blackout), despite being part of an interconnected system.

**(4) Region IV (3–5.5 EUR/kWh):** This range includes 28 of the 51 European countries considered under the baseline VOLL assumption (factor 1). The median outage rate in this group is 0.07%, below the average baseline level of 0.04% (5h of blackout), and 90% of regions in this group experience outages below 0.12%. Elevated outage levels in countries appear mostly in Eastern Europe.

**(5) Region V (5.5 EUR/kWh):** This region represents 17 countries with highly resilient electricity systems. In these cases, outages are negligible (median ≈ 0%; 80% of countries < 0.001%). Most high-income European countries fall within this category under scenarios with elevated VOLL (factors 2 or 10), reflecting a high economic valuation of supply security and strong investment signals for infrastructure.

In the European context, the Value of Lost Load is a decisive factor in determining the extent of electricity supply interruptions. A sufficiently high VOLL leads to nearly uninterrupted supply, whereas low values result in significant or even complete outages. The analysis underscores the importance of carefully choosing VOLL values in modeling studies, particularly when analyzing cross-border electricity trade and system integration.